\pdfoutput=1
\documentclass{scrartcl}
\usepackage{amsthm,amssymb,amsmath,bm,graphicx,graphbox,
  dsfont,xcolor,float,booktabs,hyperref,units,mathtools,subcaption,siunitx}
\usepackage{mycommands}
\usepackage{setspace}
\usepackage[style=numeric-comp]{biblatex}
\usepackage[Algorithm]{algorithm} %
\usepackage{algpseudocode} %
\usepackage[inline]{enumitem} %

\usepackage{tikz}
\usetikzlibrary{shapes, positioning, fit, quotes, angles, fit, calc, patterns}
\usepackage{pgfplots}
\usepgfplotslibrary{groupplots}
\pgfplotsset{compat=newest,
  select coords between index/.style 2 args={ 
    x filter/.code={
      \ifnum\coordindex<#1\fi
      \ifnum\coordindex>#2\fi
    }
  },
}

\definecolorset{RGB}{tol}{}{blue,68,119,170;cyan,102,204,238;green,34,136,51;yellow,204,187,68;red,238,102,119;purple,170,51,119;grey,187,187,187}
\pgfplotscreateplotcyclelist{bright}{%
	tolblue,every mark/.append style={fill=tolblue!80!black},mark=*\\%
	tolcyan,every mark/.append style={fill=tolcyan!80!black},mark=square*\\%
	tolgreen,every mark/.append style={fill=tolgreen!80!black},mark=otimes*\\%
	tolyellow,every mark/.append style={fill=tolyellow!80!black},mark=star\\%
	tolred,every mark/.append style={fill=tolred!80!black},mark=diamond*\\%
	tolpurple,every mark/.append style={solid,fill=tolpurple!80!black},mark=*\\%
	tolgrey,every mark/.append style={solid,fill=tolgrey!80!black},mark=square*\\%
}

\pgfplotsset{
	cycle list name=bright
}

\newcommand{\blue}[1]{{\color{blue} \noindent #1}}
\title{Deep neural networks for geometric multigrid methods} %

\author{Nils Margenberg\footnotemark[2]\and
  Robert Jendersie\footnotemark[3]\and
  Thomas Richter\footnotemark[4] \and
  Christian Lessig\footnotemark[3]}

\begin{document}
\maketitle

\begin{abstract}
  We investigate scaling and efficiency of the deep neural network multigrid method (DNN-MG).
  DNN-MG is a novel neural network-based technique for the simulation of the Navier-Stokes equations that combines an adaptive geometric multigrid solver, i.e. a highly efficient
  classical solution scheme, with a recurrent neural network with memory.
  The neural network replaces in DNN-MG one or multiple finest multigrid layers and provides a correction for the classical solve in the next time step.
  This leads to little degradation in the solution quality while substantially reducing the overall computational costs.
  At the same time, the use of the multigrid solver at the coarse scales allows for a compact network that is easy to train, generalizes well, and allows for the incorporation of physical constraints.
  Previous work on DNN-MG focused on the overall scheme and how to enforce divergence freedom in the solution.
  In this work, we investigate how the network size affects training and solution quality and the overall runtime of the computations.
  Our results demonstrate that larger networks are able to capture the
  flow behavior better while requiring only little additional training time.
  At runtime, the use of the neural network correction can even reduce the computation time compared to a classical multigrid simulation through a faster convergence of the nonlinear solve that is required at every time step.
\end{abstract}

\renewcommand{\thefootnote}{\fnsymbol{footnote}}
\footnotetext[2]{Helmut Schmidt University,
Holstenhofweg 85, 22043 Hamburg,
Germany,
\texttt{margenbn@hsu-hh.de}}
\footnotetext[3]{University of Magdeburg,
Institute for Simulation and Graphics,
Universit\"atsplatz 2,
39104 Magdeburg,
Germany,
\texttt{christian.lessig@ovgu.de}}
\footnotetext[4]{University of Magdeburg,
Institute for Analysis and Numerics,
Universit\"atsplatz 2,
39104 Magdeburg,
Germany,
\texttt{thomas.richter@ovgu.de}}

\section{Introduction}
\label{sec:org1baae5d}
In the last decade, deep neural networks had great success with tasks such as machine
translation and image classification and recently also showed surprising effectiveness for even more challenging problems such as image generation and playing games, cf.~\cite{LeCun2015}.
Underlying these results is the ability of deep neural networks (DNNs) to accurately approximate
high dimensional mappings when sufficiently deep and complex network architectures are used, training is performed on large or very large amounts of data, and powerful hardware, such as GPUs, are employed for this.

The aforementioned success of deep neural networks and related techniques leads to a growing interest
to apply these also to problems in computational science and engineering, including for the simulation of partial differential equations (PDEs).
Raissi, Karniadakis and co-workers~\cite{Raissi2018,Raissi2019b}, for example, proposed physics-informed neural networks (PINNs) that directly learn the mapping from the input of a potentially parametric PDE  to its solution.
Kasim et al.~\cite{Kasim2020} showed the effectiveness of deep neural networks to approximate a wide range of partial differential equations when also the network architecture is part of the training.
Various works~\cite{Li2020b,Lu2020,Li2020} recently also developed frameworks for operator learning. These can be applied to the fundamental solution of a PDE.
A natural connection between neural networks and dynamical systems has been observed for example by E~\cite{E2017} and Haber and Rothutto~\cite{Haber2017,Chang2018b}. This can also be exploited when neural networks are used for the solution of PDEs, for example to improve the stability of the learning.
A more detailed discussion of related work can, for example, be found in~\cite{hartmannNeuralNetworkMultigrid2020}.

A principle question for the applicability of deep neural networks for PDEs is how the equations themselves and known results on their solutions should be incorporated.
Furthermore, there exists already a large number of numerical schemes for partial differential equations and one has to ask under what circumstances a neural network can outperform these.
To address these questions, we investigated in our previous work how an adaptive multigrid solver, i.e. one of the most efficient classical techniques for partial differential equations, can be combined with deep neural networks for the solution of the Navier-Stokes equations.
The resulting deep neural network multigrid method (DNN-MG)~\cite{hartmannNeuralNetworkMultigrid2020} combines a geometric multigrid solver with a recurrent
neural network with memory to replace computations on one or multiple finest mesh levels
with network based corrections. Through this, DNN-MG achieves a significant
speed-up with only a small degradation in quality while its design for small patches ensures that the technique generalizes well and is comparatively easy to train.
The divergence freedom of learned corrections for the Navier-Stokes equations was investigated in \cite{margenbergStructurePreservationDeep2020}.

In this work, we investigate the scalability with respect to the size
of the neural network.
Our results show that larger networks are able to capture the
flow behavior better and improve important flow quantities.
The increased network size has a minor impact on training and test
time. DNN-MG compared to a pure numerical
simulation can even reduce the simulation time.

\section{Deep neural network multigrid solver}
\label{sec:orgc64f49c}
In this section we give a short summary of the deep neural network
multigrid solver~\cite{hartmannNeuralNetworkMultigrid2020}.
Although it is applicable to other equations, we will consider it in the
context of the incompressible Navier-Stokes equations as in our previous work.

\subsection{Finite Element Discretization of the Incompressible Navier-Stokes equations}
\label{sec:org1a8d80c}
The incompressible Navier-Stokes equations are given by
\begin{equation}
  \begin{alignedat}{2}
    \label{eq:nsstrong}
    \partial_t v + (v\cdot \nabla)v - \frac{1}{\mathrm{Re}}\Delta v
    +\nabla p &= f \quad &&\text{on } [0,\,T] \times \Omega\\
    \nabla \cdot v &= 0 \quad &&\text{on } [0,\,T] \times \Omega.
  \end{alignedat}
\end{equation}
Here \(v\colon [0,\,T]\times \Omega \to \R^2\) is the velocity,
\(p\colon [0,\,T]\times \Omega \to \R\) the pressure,
\(\mathrm{Re}>0\) the Reynolds number, and \(f\) an external force.
The initial and boundary conditions are
\begin{equation}
  \label{eq:boundary}
  \begin{alignedat}{2}
    v(0,\,\cdot) &= v_0(\cdot)\quad &&\text{on }\Omega\\
    v &= v^D \quad &&\text{on } [0,\,T] \times \Gamma^D\\
    \frac{1}{\mathrm{Re}}(\vec n\cdot\nabla)v - p\vec n &=0 \quad &&\text{in } [0,\,T]
    \times \Gamma^N
  \end{alignedat}
\end{equation}
where \(\vec n\) is the outward facing unit normal on the boundary \(\partial\Omega\)
of the domain. On the outflow boundary \(\Gamma^N\) we consider the do-nothing outflow
condition~\cite{heywoodArtificialBoundariesFlux1992} which is a well established model
for artificial boundaries; see~\cite{braackDirectionalDonothingCondition2014}
for a detailed discussion.

We use a weak finite element formulation to discretize~\eqref{eq:nsstrong} with
\(\smash{v_h,\,\phi_h \in V_h = [W_h^{(2)}]^d}\)
and test functions \(\smash{p_h,\,\xi_h \in L_h = W_h^{(2)}}\).
\(\smash{W_h^{(r)}}\) is here the space of continuous functions which
are polynomials of degree \(r\) on each element
\(T \in\Omega_h\) and \(\Omega_h\) is the mesh domain.
The equal order finite element pair \(V_h\times L_h\) that we use does not fulfill the
inf-sup condition. We hence use stabilization terms of local projection
type~\cite{beckerFiniteElementPressure2001} with parameter \(\alpha_T = \alpha_0\cdot  \mathrm{Re}
   \cdot h_T^2\) and projection \(\smash{\pi_h:W_h^{(2)}\to W_h^{(1)}}\) into the space
of linear polynomials.

We use the second order Crank-Nicolson method for time
discretization.
The solution of~\eqref{eq:nsstrong} subject to~\eqref{eq:boundary}
at the $n$-th time step then amounts to determining the state
$x_n=(v_n^{1},\dots,\,v_{n}^{d},\,p_n)$ such that
\begin{equation}
  \label{eq:nonlinearshort}
  {\cal A}_n(x_n) = f_n
\end{equation}
\begin{equation}\label{varform}
  \begin{aligned}
    [{\cal A}_n(x)]_i &\coloneqq (\nabla \cdot v_{n}, \, \xi_h^i)+
    \sum_{T\in\Omega_h}\alpha_T (\nabla (p_n-\pi_h p_n),\nabla (\xi_h^i-\pi_h \xi_h^i))
    \\[0pt]
      [f_n]_i &\coloneqq 0
      \\
  \end{aligned}
\end{equation}
for all pressure unknowns $i=1,\dots,N_p$ and
\begin{equation}\label{varformv}
  \begin{aligned}
    [{\cal A}_n(x_{n})]_{i+N_p} &\coloneqq
    \frac{1}{k}(v_n,\,\phi_h^i)\,
    +{}\frac{1}{2} (v_n\cdot \nabla v_n,\,\phi_h^i)
    +\frac{1}{2\mathrm{Re}}(\nabla v_n,\,\nabla \phi_h^i)
    -(p_n,\,\nabla \cdot \phi_h^i)
    \\[6pt]
    [f_n]_{i+N_p} &\coloneqq  \frac{1}{k}(v_{n-1},\, \phi_h^i)
    +\frac{1}{2}(f_n,\phi_h^i)
    +\frac{1}{2}(f_{n-1},\,\phi_h^i)
    \\
    &\qquad
    -\frac{1}{2}(v_{n-1}\cdot\nabla v_{n-1},\,\phi_h^i)
    - \frac{1}{2\mathrm{Re}}(\nabla v_{n-1},\,\nabla \phi_h^i),
  \end{aligned}
\end{equation}
for all velocity unknowns $i=1,\dots,N_v$ with corresponding test
functions $\phi_h^i$ and $\xi_h^i$.
Equation~\eqref{eq:nonlinearshort} is a large nonlinear system of algebraic equations that is
solved by Newton's method with the initial guess \(x_n^{(0)}=(v_{n-1},p_{n-1})\).
Each Newton step takes the form
\begin{equation}
  \label{eq:newton}
  {\cal A}_n'(x_n^{(l-1)}) \, w^{(l)} = f_n-{\cal A}_n(x_n^{(l-1)}), \quad
  x_n^{(l)}=x_n^{(l-1)}+w^{(l)} \text{ for } l=1,2,\dots
\end{equation}
\({\cal A}'(x^{(l-1)})\) is the Jacobian of \({\cal A}\) at \(x^{(l-1)}\),
which can be computed analytically for this problem,
cf.~\cite[Sec. 4.4.2]{richterFluidstructureInteractionsModels2017}.

\subsection{Deep Neural Network Multigrid Solver}
\label{sec:orga419e94}

For the solution of the Newton iteration described above we use GMRES and the geometric multigrid
method as preconditioner, leading to a highly efficient solver. However,
multiple GMRES steps and one up- and down-sweep of the multigrid method per Newton
step still result in a significant amount of computations.
Most of the computation time is thereby spent on the finest mesh level (about 75\%), so solving on
an additional level \(L+1\) leads to a substantial increase in the computational costs.
The principle idea of the deep neural network multigrid solver (DNN-MG) is therefore to use the multigrid hierarchy and correct a classical, coarse solution of the nonlinear Navier-Stokes equation on one or multiple finest mesh levels with a neural network, see Algorithm~\ref{alg:dnnmg}.
For this to be meaningful, the network has to be able to (largely) retain the accuracy of the multigrid computation.
Through the coarse multigrid solution as well as the nonlinear residual, the network has, however, highly informative input; in other words, the classical coarse solution as well as its error on finer levels provide a strong prior for the network correction.
We furthermore use the underlying mesh to subdivide our domain into patches, e.g. a mesh element or a structured assembly of a few adjacent elements, and apply the network on these.
This greatly aids generalizability since the network no longer needs to correct the general flow and instead only local corrections for each patch are required.
Furthermore, the local setup simplifies training and ensures that a large training corpus is given by just a few example flows.
As neural network we use a recurrent one with memory, more specifically GRUs, so that complex flow behavior can be predicted and coherence in time is ensured.
%
%
\begin{algorithm}[t]
  \setstretch{1.2}
  \caption{DNN-MG for the solution of the Navier-Stokes equations.
    Lines 6--9 (blue) provide the modifications of the DNN-MG method compared
    to a classical Newton-Krylow simulation with geometric multigrid preconditioning.}
  \label{alg:dnnmg}
  \begin{algorithmic}[1]
    \For{all time steps $n$}
    \While{not converged}\Comment{Newton's method for Eq.~\ref{varform}}
      \State{$\delta z_i \leftarrow$ \Call{multigrid}{$L,\,A_{L}^n,\,b_{L}^n,\,\delta
          z_i$}}\Comment{Geometric multigrid with $z = (p_n^L,\,v_n^L)$}
      \State{$z_{i+1} \leftarrow z_i + \epsilon \, \delta z_i$}
    \EndWhile{}
    \State{\blue{$\tilde{v}_n^{\scriptscriptstyle L+1} \leftarrow \mathcal{P}(v_n^{\scriptscriptstyle L}) $}}\Comment{\blue{Prolongation on level $L+1$}}
    \State{\blue{$d_n^{\scriptscriptstyle L+1} \leftarrow \mathcal{N}(\tilde{v}_n^{\scriptscriptstyle L+1},\,\Omega_L,\,\Omega_{L+1})$}}\Comment{\blue{Prediction of velocity correction}}
    \State{\blue{$b_{n+1}^{\scriptscriptstyle L+1} \leftarrow \textbf{Rhs}(\tilde{v}_n^{\scriptscriptstyle L+1} + d_n^{\scriptscriptstyle L+1},f_n,f_{n+1})$}}\Comment{\blue{Set up rhs for next time step}}
    \State{\blue{$b_{n+1}^{\scriptscriptstyle L} \leftarrow \mathcal{R}(b_{n+1}^{\scriptscriptstyle L+1})$}}\Comment{\blue{Restriction of rhs to level $L$}}
    \EndFor{}
  \end{algorithmic}
\end{algorithm}
In summary, the key features that make DNN-MGs efficient and flexible are:
\begin{enumerate}[label=(\alph*)]
\item patchwise operation to ensures generalizability to unseen flow regimes and meshes\label{item:nn:i};
\item GRUs with memory to capture complex flow behavior and ensure coherence in time\label{item:nn:ii};
\item use of the nonlinear residual in Eq.~\ref{eq:nonlinearshort} as network input to have rich information
about the sought correction.\label{item:nn:iii}
\end{enumerate}

Algorithmically, DNN-MG works as follows, see Alg.~\ref{alg:dnnmg}.
The neural network-based correction is applied at the end of every time step after we computed an
updated velocity \(v_n^L\) on level \(L\) (Alg.~\ref{alg:dnnmg}, l. 2-4).
For this, \(v_n^L\) is first prolongated to level \(L+1\),
yielding \(\tilde{v}_n^{L+1} \coloneqq \mathcal{P}(v_n^L)\).
Then we compute the input to the neural network, which is evaluated individually for each patch \(P_i\)
(a mesh element on level \(L\) in the present work).
The inputs to the neural network are thereby also entirely local and include
the residual of Eq.~\eqref{eq:nonlinearshort} on level \(L+1\), the prolongated velocity
\(\smash{\tilde{v}_n^{L+1}}\), geometric properties such as the patch's aspect ratio,
and the P\'eclet number over the patch.
The network predicts a velocity correction \(d_{n,i}^{L+1}\) and notably does not include the
pressure. It enters only through the fine mesh residual, which is part of the input.

At the end of the neural network-based correction a provisional right hand side
\(b_{n+1}^{\scriptscriptstyle L+1} = \mathrm{Rhs}(\tilde{v}_n^{\scriptscriptstyle L+1} +
   d_n^{\scriptscriptstyle L+1},f_n,f_{n+1})\)
of~\eqref{eq:newton} is computed on level \(L+1\) and then restricted to level \(L\).
The corrected \(b_{n+1}^{\scriptscriptstyle L}\) is then used in the
next time step, which is again improved by a neural network-based
correction at the end of the step.
Through this, the neural network based correction propagates back into the
Newton solver and improves the numerical solution.

Before we can use DNN-MG, we train the network with a high fidelity finite element
solution on the fine mesh level \(L+1\) using the loss function
\begin{equation}\label{loss}
  \mathcal{L}(v^L,v^{L+1};d^{L+1})
  \coloneqq \sum_{n=1}^{N}\sum_{P_i\in \Omega_{L+1}}
  \big\Vert v^{L+1}_{n+1} - (v^L_n+d^{L+1}_n)\big\Vert^2_{l^2(P_i)}.
\end{equation}
Here \(N\) is the number of time-steps in the training data.
Since DNN-MG operates strictly local on patches, the loss accumulates
the local residuals for each patch in the second sum.

For a detailed description and results showing the generalizability, efficiency and
efficacy, we refer to
\cite{hartmannNeuralNetworkMultigrid2020,margenbergStructurePreservationDeep2020}.

\section{Numerical experiments}
\label{sec:orgc3dd7d1} 

In this section we investigate the scalability of DNN-MG with respect to the
size of the neural network both in terms of its computational efficiency
and the quality of the results.

\subsection{Setup}

As test problem we considered the classic
benchmark of a laminar flow around an
obstacle~\cite{SchaeferTurek1996} governed by the Navier-Stokes equations.
Analogous to~\cite{hartmannNeuralNetworkMultigrid2020} and
different from~\cite{SchaeferTurek1996}
we use elliptic obstacles with varying aspect ratios instead of
discs. This avoids that the neural networks can memories flows
and perform well in the experiments.

The neural networks of DNN-MG used for our experiments had the same overall structure
as in our previous work~\cite{hartmannNeuralNetworkMultigrid2020} with GRU cells
with memory at their heart.
In contrast to~\cite{hartmannNeuralNetworkMultigrid2020}, however, we replaced the
convolutional layer with a dense one since it slowed
down training considerably with only a negligible effect on the effectiveness of the
network.
To investigate the effect of the neural network size we
parameterized the network by the hidden state size of the GRU-cells (32, 64) and the
number of GRUs stacked on top of each other (1, 2, 3).
The resulting network sizes are reported in Tab.~\ref{tab:training}.

We use the finite element
library Gascoigne 3d~\cite{beckerFiniteElementToolkit} for numerical simulations
and PyTorch~\cite{paszkeAutomaticDifferentiationPyTorch2017} for the neural networks.
The training data was obtained in the same manner described in
previous works by running highly resolved finite element simulations
consisting of 2560 elements and 8088 degrees of freedom
and localizing the data to the time steps and to the patches of the
mesh.
The simulations were carried out on a coarser mesh with 640 elements
resulting 2124 degrees of freedom.
%
The training of the networks was performed on a NVIDIA Tesla V100.
The numerical simulations were performed on the same machine equipped with
2 Intel Xeon E5-2640 v4 CPUs and running on 20 cores.

\subsection{Training performance}

\begin{figure}[t]
  \includegraphics[width=\textwidth]{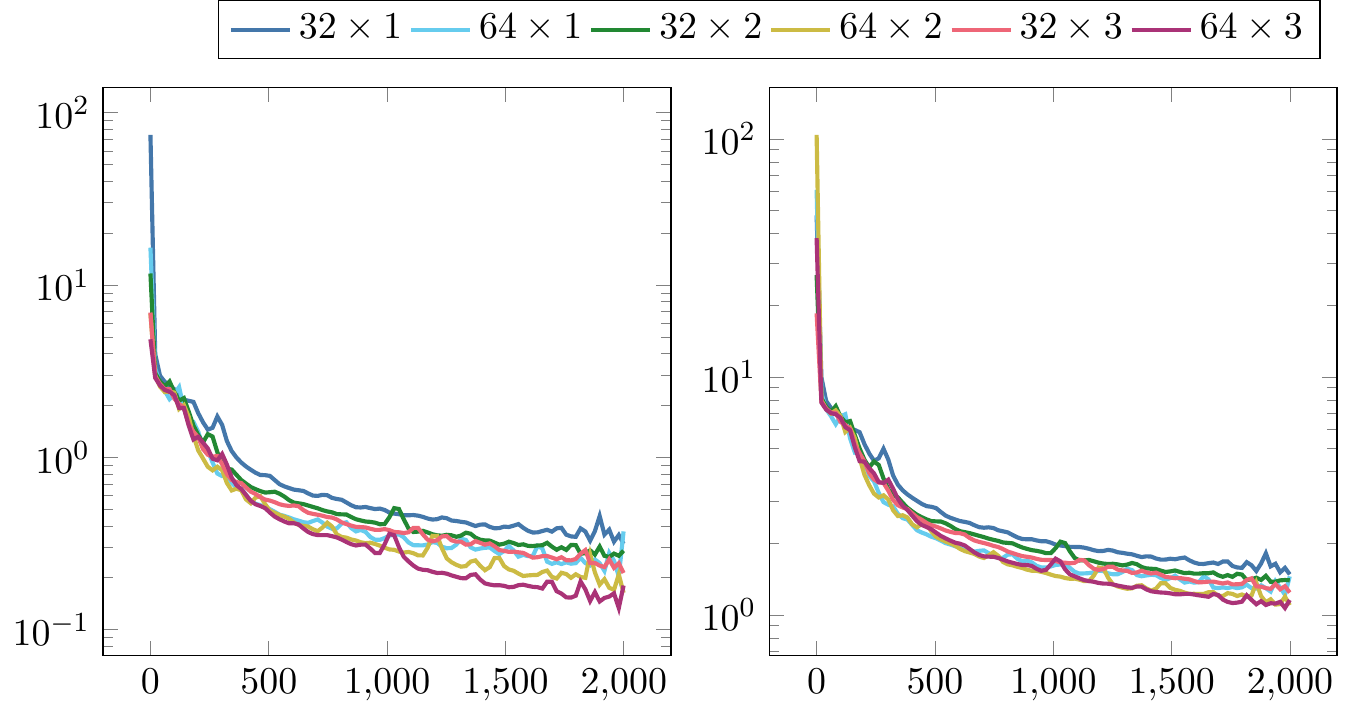}
  \caption{Average training (left) and validation (right) loss of different
    GRU-cells per patch (smoothed with gnuplots smooth bezier curve): A
    cell with dimension $m\times n$ consists of $n$ GRUs stacked on top of each
    other, where the hidden states $h$ are of size $m$.
  }\label{fig:losses}
\end{figure}

The convergence behavior of the training for the different network
configurations is reported in Fig.~\ref{fig:losses}.
We see that all networks show a similar behavior although the
larger networks converge slightly faster and the
loss is generally lower. The strongest effect on this behavior has the
size of the hidden state with the GRUs of size $64\times n$ performing best.

In Tab.~\ref{tab:training} we see that the training time scales very well with the
network size. Although the number of parameters is increased
by almost one magnitude between the networks $32 \times 1$ and $64 \times 3$, the training time grows only by
35\%.

\begin{table}[ht]
  \parbox[c]{.49\linewidth}{ \centering
    \begin{tabular}{ccccc}
      \toprule
      GRU dims & \multicolumn{2}{c}{\#parameters} & \multicolumn{2}{c}{time}\\
      \midrule
      $32\times 1$ & 8544  & $(\num{1.00})$ & \SI{7039}{\second} & $(\num{1.00})$\\
      $64\times 1$ & 23232 & $(\num{2.72})$ & \SI{7238}{\second} & $(\num{1.02})$\\
      $32\times 2$ & 14688 & $(\num{1.72})$ & \SI{8055}{\second} & $(\num{1.14})$\\
      $64\times 2$ & 47808 & $(\num{5.60})$ & \SI{8335}{\second} & $(\num{1.18})$\\
      $32\times 3$ & 20832 & $(\num{2.44})$ & \SI{9158}{\second} & $(\num{1.30})$\\
      $64\times 3$ & 72384 & $(\num{8.47})$ & \SI{9573}{\second} & $(\num{1.36})$\\
      \bottomrule
    \end{tabular}}
  \parbox[c]{.49\linewidth}{ \centering
    \begin{tikzpicture}
      \begin{axis}[xlabel=\#parameters ($\times 10^4$),ylabel={time [\,s\,]},%
        width=0.45\textwidth, ymin=0,ymax=10000,
        tick scale binop=\times,
        every axis plot/.append style={very thick},%
        legend style={at={(0.5,0.1)},anchor=south}]
        \addplot coordinates{ %
          (.8544, 7039) %
          (1.4688, 8055) %
          (2.0832, 9158) %
        }; %
        \addplot coordinates{ %
          (2.3232, 7238) %
          (4.7808, 8335) %
          (7.2384, 9573) %
        };
        \legend{$32\times n$, $64\times n$}
      \end{axis}
    \end{tikzpicture}}
  \caption{GRU dimensions, parameters and the training time. In parentheses we list
    the factor compared to GRU $32\times 1$.}\label{tab:training}
\end{table}

\subsection{Application performance}

Once trained, the network is used in the DNN-MG solver described
in Section~\ref{sec:orgc64f49c}.
Table~\ref{tab:timing} collects the runtimes for the overall
simulation over 1050 time steps as well as the time required for
evaluating the neural network in DNN-MG. In
most cases, the neural network contributes less than 4\% to the
overall computational time. We further note that about
$\SI{1.5}{\second}$ of the network evaluation time must be attributed
to data processing and communication between the finite element library
and PyTorch.
In an optimized implementation this time could most likely be reduced
substantially.

The observed scaling in Tab.~\ref{tab:timing} does not always agree with what one
would expect from the neural networks sizes in Tab.~\ref{tab:training}.
There are multiple reasons for this.
First, when the neural network correction is highly inaccurate then we need more Newton
iterations in the next time step, resulting in higher runtimes such as for the
$32\times 1$ network. High quality predictions, on the other hand, improve the runtimes
and for the $64\times 3$ GRU-cell, for example, the time associated with the
assembly routines was reduced by a factor of 4.
Second, if the network deteriorates the solution and convergence of the Newton iteration is too slow
we switch off the prediction for this time step. If this happens often, it leads
to significantly less network evaluations and thus reduces runtimes of the
network (this causes the spikes e.g. for $64\times 1$, $64\times 2$)
although at the price of a loss in accuracy.

In summary, the networks take up only a small fraction of the runtime and may
even reduce the overall runtime through faster convergence of the Newton
iteration.
For the largest network $64 \times 3$, for example, we achieve a speedup of almost 100\%
compared to the solution without a network.
However, an ill-suited network can also negatively affect the convergence and hence
slow down computations.
\begin{table}[h]
  \parbox[c]{.49\linewidth}{ \centering
    \begin{tabular}{cccc}
      \toprule
      type & runtime & \multicolumn{2}{c}{NN-Eval}\\
      \midrule
      MG($L+1$)         & \SI{538.45}{\second} &  & $(\num{0}\%)$\\
      MG($L$)           & \SI{453.06}{\second} &  & $(\num{0}\%)$\\
      $32\times 1$ & \SI{705.53}{\second} & \SI{21.84}{\second} & $(\num{3.1}\%)$\\
      $64\times 1$ & \SI{284.63}{\second} & \SI{2.85}{\second}  & $(\num{1.0}\%)$\\
      $32\times 2$ & \SI{257.72}{\second} & \SI{7.44}{\second}  & $(\num{2.9}\%)$\\
      $64\times 2$ & \SI{224.95}{\second} & \SI{4.14}{\second}  & $(\num{1.8}\%)$\\
      $32\times 3$ & \SI{200.42}{\second} & \SI{11.01}{\second} & $(\num{5.5}\%)$\\
      $64\times 3$ & \SI{230.6}{\second} & \SI{5.47}{\second}   & $(\num{2.3}\%)$\\
      \bottomrule
    \end{tabular}}
  \parbox[c]{.49\linewidth}{ \centering
    \begin{tikzpicture}
      \begin{axis}[xlabel=\#parameters ($\times 10^4$),ylabel={time [\,s\,]},%
        width=0.45\textwidth,%
        every axis plot/.append style={very thick},%
        legend style={at={(0.75,0.5)},anchor=south}]
        \addplot coordinates{ %
          (.8544,  21.84) %
          (1.5120, 7.44) %
          (2.0832, 11.01) %
        }; %
        %
        \addplot coordinates{ %
          (2.3232, 2.85) %
          (4.7808, 4.14) %
          (7.2384, 5.47) %
        }; \legend{$32\times n$, %
          $64\times n$}
      \end{axis}
    \end{tikzpicture}}
  \caption{Timing results for the simulations with the different GRUs we tested with
   \emph{runtime} referring to the timing of the whole simulation and \emph{NN-eval}
    the time spent on evaluating the neural network. In general $\approx
    \SI{1.5}{\second}$ of the latter is spent on data processing. In parentheses we
    specify the time taken relative to the runtime in percent.
  }\label{tab:timing}
\end{table}

\subsection{Accuracy of the networks}

To compare the approximation power of the different neural
network configurations, we evaluate the solution to the
benchmark problem with respect to three functionals and also
compute the error in the computed velocities and pressures.

The first functional we consider is the squared divergence
\[
J_{\text{div}}(v(t)) 
= \int_\Omega |\nabla\cdot v(t)|^2\,\text{d}x.
\]
We know that the exact solution should satisfy
$J_{\text{div}}(v)=0$. However, since our finite element approach is not
strictly divergence free and since we are adding artificial
stabilization terms, compare~(\ref{varform}), the numerical solution
will carry a divergence error that is also present in the training data.
Hence, also the neural network correction will not be exactly
divergence free although it should not deteriorate it further.

As suggested in~\cite{SchaeferTurek1996}, the second and third functional
we are considering are the drag and lift forces acting on the obstacle
\[
\begin{aligned}
J_{\text{drag}}(v(t),p(t))& = -\int_{\partial\Gamma} \big(\frac{1}{Re}\nabla v(t) -
p(t)I\big)\vec n\cdot \vec e_x\,\text{d}s,\\
J_{\text{lift}}(v(t),p(t))& = -\int_{\partial\Gamma} \big(\frac{1}{Re}\nabla v(t) -
p(t)I\big)\vec n\cdot \vec e_y\,\text{d}s,
\end{aligned}
\]
where we denote by $\Gamma$ the boundary of the obstacle, by $\vec n$
the unit normal vector pointing into the obstacle and by $\vec
e_x,\vec e_y$ the Cartesian unit vectors.

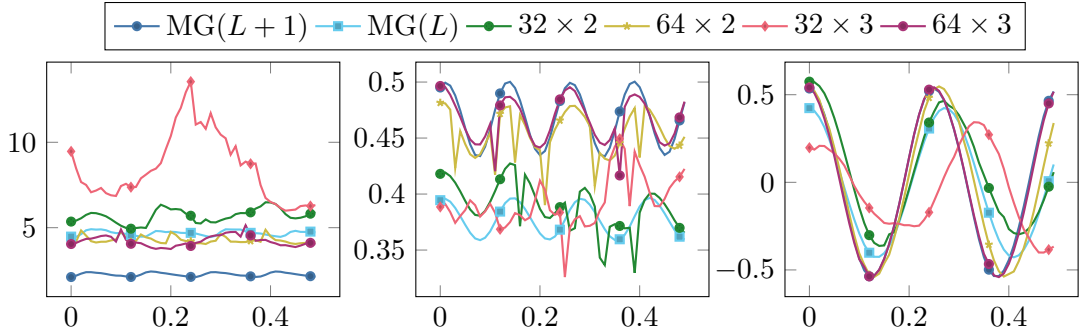
\begin{figure}[ht]
\begin{tikzpicture} [every mark/.append style={mark size=1.5pt}]
\begin{groupplot}[group style={columns=3}, width=0.37\textwidth,
	legend columns=-1,
	legend style={at={(-0.75,1.05)},anchor=south},
  every axis plot/.append style={thick}]
	\pgfplotsinvokeforeach {1, 2 ,3} {
			\nextgroupplot
			\addplot+[ mark repeat=12] table[x index = 0, y index = #1, col sep=space ] {data/fcts_ctest_allMGL1.txt};
			\addplot+[ mark repeat=12] table[x index = 0, y index = #1, col sep=space] {data/fctsctestallMG.txt};
      \addplot+[ mark repeat=12] table[x index = 0, y index = #1, col sep=space ] {data/fctsctestallDNNMG32times2.txt};
      \addplot+[ mark repeat=12] table[x index = 0, y index = #1, col sep=space ] {data/fctsctestallDNNMG64times2.txt};
      \addplot+[ mark repeat=12] table[x index = 0, y index = #1, col sep=space ] {data/fctsctestallDNNMG32times3.txt};
      \addplot+[ mark repeat=12] table[x index = 0, y index = #1, col sep=space] {data/fctsctestallDNNMG64times3.txt};
    }
    \legend{MG($L+1$),%
      MG($L$),%
      $32\times 2$,%
      $64\times 2$,%
      $32\times 3$,%
      $64\times 3$%
    }
\end{groupplot}
\end{tikzpicture}
\caption{Divergence $J_{\text{div}}(v(t))$, drag $J_{\text{drag}}(v(t),p(t))$, lift $J_{\text{lift}}(v(t),p(t))$ (from left to right) for the coarse multigrid solution MG($L$), a high fidelity
  finite element reference solution MG($L+1$) and DNN-MG solutions
  for different GRU configurations.}
\label{fig:func}
\end{figure}
Figure~\ref{fig:func} shows all three functionals over the temporal
interval $[\unit[9]{s},\unit[10]{s}]$, where the solution reached the
periodic limit cycle (the phase were adjusted so that they agree
on the first maximum in $[\unit[9]{s},\unit[10]{s}]$).
We use a high resolution
simulation MG($L+1$) as reference and MG($L$), that is the uncorrected
simulation on level $L$, as base line.
We observe that most functional outputs are
significantly improved by DNN-MG
with the best performing neural networks.
Nearly exact values are obtained in case of the lift functional (right) but also drag (middle) is
getting close to the reference value, in particular for $64 \times 3$.
For the divergence $J_{\text{div}}$ the improvement is more
modest although compared to our previous results~~\cite{ETNAmargenbergStructurePreservationDeep2020}
(\cite[Fig. 1]{margenbergStructurePreservationDeep2020} in the
corresponding preprint) with a network with $32\times 1$ GRU cells we
still obtain a lower divergence.

\begin{figure}[t]
\begin{subfigure}[t]{.495\textwidth}
	\centering
	\begin{tikzpicture} [every mark/.append style={mark size=1.5pt}]
		\begin{axis}[xlabel={time $t$},
      ylabel=velocity error, width=1.0\textwidth,
      legend pos=north west,
      every axis plot/.append style={thick}]
      \addplot+[ mark repeat=12] table[x index = 0, y index = 2, col sep=space ] {data/vperrMG.txt}; %
      \addplot+[ mark repeat=12] table[x index = 0, y index = 2, col sep=space ] {data/vperrDNNMG32times2.txt}; %
      \addplot+[ mark repeat=12] table[x index = 0, y index = 2, col sep=space ] {data/vperrDNNMG64times2.txt}; %
      \addplot+[ mark repeat=12] table[x index = 0, y index = 2, col sep=space ] {data/vperrDNNMG32times3.txt}; %
      \addplot+[ mark repeat=12] table[x index = 0, y index = 2, col sep=space ] {data/vperrDNNMG64times3.txt}; %
      \legend{MG,
        $32\times 2$,%
        $64\times 2$,%
        $32\times 3$,%
        $64\times 3$%
      }
		\end{axis}
	\end{tikzpicture}
\end{subfigure}
\begin{subfigure}[t]{.495\textwidth}
	\centering
	\begin{tikzpicture} [every mark/.append style={mark size=1.5pt}]
		\begin{axis}[xlabel={time $t$},
      ylabel=pressure error, width=1.0\textwidth,
      legend pos=north west,
      every axis plot/.append style={thick}]
      \addplot+[ mark repeat=12] table[x index = 0, y index = 2, col sep=space ] {data/vperrMG.txt}; %
      \addplot+[ mark repeat=12] table[x index = 0, y index = 2, col sep=space ] {data/vperrDNNMG32times2.txt}; %
      \addplot+[ mark repeat=12] table[x index = 0, y index = 2, col sep=space ] {data/vperrDNNMG64times2.txt}; %
      \addplot+[ mark repeat=12] table[x index = 0, y index = 2, col sep=space ] {data/vperrDNNMG32times3.txt}; %
      \addplot+[ mark repeat=12] table[x index = 0, y index = 2, col sep=space ] {data/vperrDNNMG64times3.txt}; %
      \legend{MG,
        $32\times 2$,%
        $64\times 2$,%
        $32\times 3$,%
        $64\times 3$%
      }
		\end{axis}
	\end{tikzpicture}
\end{subfigure}
\caption{Relative errors for velocity (left) and pressure (right) compared to the reference solution on level $L + 1$ for the coarse solution on level $L$ and DNN-MG with different neural network configurations.}
\label{fig:norms}
\end{figure}
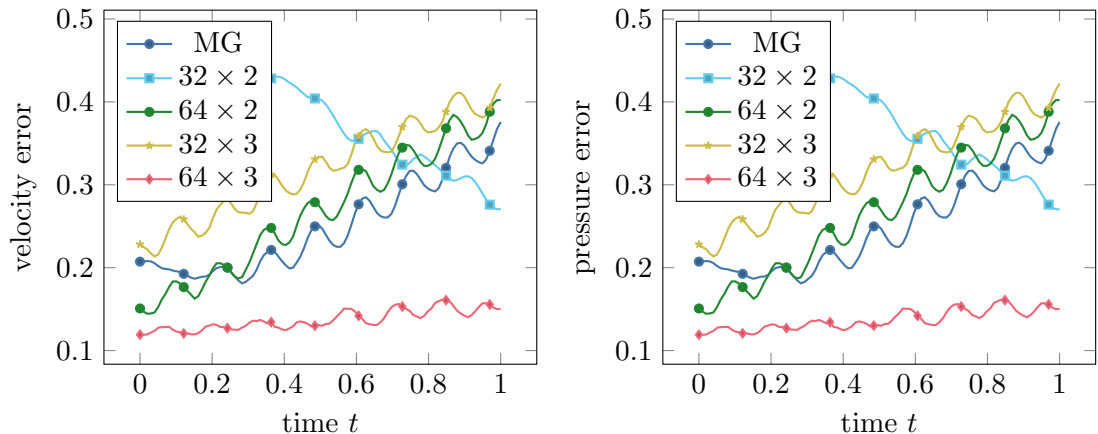

In Fig.~\ref{fig:norms} we show the relative $l_2$ errors of the
velocity and pressure with respect to the fine mesh solution MG($L+1$).
The results demonstrate that the largest network $64 \times 3$
is able to produce robust and accurate predictions for large time intervals.
This network is also able to substantially reduce the phase error that
results from the interplay between spatial and temporal
discretization~\cite{richterParallelTimesteppingFluidstructure2021,hartmannNeuralNetworkMultigrid2020}
and leads to lower frequency oscillations on coarser meshes.
Our previous results and the smaller networks show the same effect for the DNN-MG
approach (low frequency oscillations in Fig.~\ref{fig:norms}) but using a larger neural network cures this
defect. We also refer to Fig.~\ref{fig:func} which shows (e.g. in the
case of the lift functional) that the GRU $64\times 3$ solution is almost perfectly
in phase with the fine mesh solution MG($L+1$) while a lower frequency
is found for the coarse mesh MG($L$).

\section{Discussion}

The presented results show the scalability of the DNN-MG method
with respect to the network size. We observed that larger
networks lead to DNN-MG simulations that capture flow behavior more accurately.
At the same time, the increased number of parameters in the network has
little impact on the training times and can (through the faster convergence
of the Newton solve) even improve the runtime performance compared to a standard MG solver.
Smaller networks are often not able to capture the flow behavior and even render
the solution useless at times. In some cases they also have a negative impact on the
runtime performance due to the interaction of the neural network and the Newton solve.

\section{Conclusion}
\label{sec:org79b5604}
We have investigated the scalability of the DNN-MG method with respect to the
neural network size in terms of accuracy, runtime, and training time.
We showed that larger networks are able to more accurately predict flow
features and this allowed us to substantially improve over previous results for the
DNN-MG method.
In particular, our largest network $64 \times 3$ is able to recover the flow frequency of the
high fidelity simulations and through this we obtained velocity and pressure errors not attainable before.
The divergence is now also at least as
good as for the standard multigrid simulation without any additional effort, cf.~\cite{margenbergStructurePreservationDeep2020}.
Importantly, the improved accuracy comes at virtually no additional cost in terms of runtime
and are now able to improve over a standard multigrid simulation by almost 100\%.

The presented results lead to further research directions we would like to investigate in
the future.
One question is when the network size saturates and how the required amount of training
data scales with the number of network parameters.
The presented results in 2D provide also a promising basis for a hybrid finite element / neural network
simulation for the Navier-Stokes equations in 3D which we would like to consider in the future.
Finally, our practical findings provide a strong motivation to address more theoretical
questions such as the stability, consistency and convergence of DNN-MG.

\paragraph{Acknowledgement}
NM acknowledges support by the Helmholtz-Gesellschaft grant number HIDSS-0002 DASHH.

\printbibliography
\end{document}